# EXPONENTIAL INEQUALITIES FOR SELF-NORMALIZED MARTINGALES WITH APPLICATIONS

By Bernard Bercu and Abderrahmen Touati

*Université Bordeaux 1 and Faculté des Sciences de Bizerte*

We propose several exponential inequalities for self-normalized martingales similar to those established by De la Peña. The keystone is the introduction of a new notion of random variable heavy on left or right. Applications associated with linear regressions, autoregressive and branching processes are also provided.

**1. Introduction.** Let $(M_n)$ be a locally square integrable real martingale adapted to a filtration $\mathbb{F} = (\mathcal{F}_n)$ with $M_0 = 0$. The predictable quadratic variation and the total quadratic variation of $(M_n)$ are respectively given by

$$\langle M \rangle_n = \sum_{k=1}^n \mathbb{E}[\Delta M_k^2 | \mathcal{F}_{k-1}] \quad \text{and} \quad [M]_n = \sum_{k=1}^n \Delta M_k^2$$

where $\Delta M_n = M_n - M_{n-1}$. The celebrated Azuma–Hoeffding inequality [4, 16, 18] is as follows.

THEOREM 1.1 (Azuma–Hoeffding's inequality). *Let $(M_n)$ be a locally square integrable real martingale such that, for each $1 \leq k \leq n$, $a_k \leq \Delta M_k \leq b_k$ a.s. for some constants $a_k < b_k$. Then, for all $x \geq 0$,*

$$(1.1) \qquad \mathbb{P}(|M_n| \geq x) \leq 2 \exp\left(-\frac{2x^2}{\sum_{k=1}^n (b_k - a_k)^2}\right).$$

Another result which involves the predictable quadratic variation $(\langle M \rangle_n)$ is the so-called Freedman inequality [13].

THEOREM 1.2 (Freedman's inequality). *Let $(M_n)$ be a locally square integrable real martingale such that, for each $1 \leq k \leq n$, $|\Delta M_k| \leq c$ a.s. for* some constant $c > 0$. Then, for all $x, y > 0$,

$$(1.2) \qquad \mathbb{P}(M_n \geq x, \langle M \rangle_n \leq y) \leq \exp\left(-\frac{x^2}{2(y + cx)}\right).$$

Over the last decade, extensive study has been made to establish exponential inequalities for $(M_n)$ relaxing the boundedness assumption on its increments. On the one hand, under the standard Bernstein condition that for $n \geq 1$, $p \geq 2$ and for some constant $c > 0$,

$$\sum_{k=1}^{n} \mathbb{E}[|\Delta M_k|^p | \mathcal{F}_{k-1}] \leq \frac{c^{p-2} p!}{2} \langle M \rangle_n,$$

Pinelis [21] and De la Peña [8] recover (1.2). Van de Geer [12] also proves (1.2) replacing $\langle M \rangle_n$ by a suitable increasing process. On the other hand, if $(M_n)$ is conditionally symmetric which means that for $n \geq 1$, the conditional distribution of $\Delta M_n$ given $\mathcal{F}_{n-1}$ is symmetric, then De la Peña [8] establishes the nice following result.

THEOREM 1.3 (De la Peña's inequality). *Let $(M_n)$ be a locally square integrable and conditionally symmetric real martingale. Then, for all $x, y > 0$,*

$$(1.3) \qquad \mathbb{P}(M_n \geq x, [M]_n \leq y) \leq \exp\left(-\frac{x^2}{2y}\right).$$

Some extensions of the above inequalities in a more general framework including discrete-time martingales can also be found in [9, 11] where the conditionally symmetric assumption is still required for (1.3). We also refer the reader to the recent survey of De la Peña, Klass and Lai [10].

By a careful reading of [8], one can see that (1.3) is a two-sided exponential inequality. More precisely, if $(M_n)$ is conditionally symmetric, then, for all $x, y > 0$,

$$(1.4) \qquad \mathbb{P}(|M_n| \geq x, [M]_n \leq y) \leq 2 \exp\left(-\frac{x^2}{2y}\right).$$

By comparing (1.4) and (1.1), we are only halfway to Azuma–Hoeffding's inequality which holds without the total quadratic variation $[M]_n$.

The purpose of this paper is to establish several exponential inequalities in the spirit of the original work of De la Peña [8]. In Section 2, we shall propose two-sided exponential inequalities involving $\langle M \rangle_n$ as well as $[M]_n$ without any assumption on the martingale $(M_n)$. Section 3 is devoted to the introduction of a new concept of random variables heavy on left or right. This notion is really useful if one is only interested in obtaining a one-sided exponential inequality for $(M_n)$. It also provides a clearer understanding of De la Peña's conditional symmetric assumption. We shall show in Section 4 that this new concept allows us to prove (1.3). As in [8], we shall also propose exponential inequalities for $(M_n)$ self-normalized by $[M]_n$ or $\langle M \rangle_n$. Section 5 is devoted to applications on linear regressions, autoregressive and branching processes. All technical proofs are postponed to the Appendix.

**2. Two-sided exponential inequalities.** This section is devoted to two-sided exponential inequalities involving $\langle M \rangle_n$ and $[M]_n$. We start with the following basic lemma.

LEMMA 2.1. *Let $X$ be a square integrable random variable with mean zero and variance $\sigma^2 > 0$. For all $t \in \mathbb{R}$, denote*

$$(2.1) \qquad L(t) = \mathbb{E}\left[\exp\left(tX - \frac{t^2}{2}X^2\right)\right].$$

*Then, we have for all $t \in \mathbb{R}$,*

$$(2.2) \qquad L(t) \leq 1 + \frac{t^2}{2}\sigma^2.$$

PROOF. The proof is given in Appendix A.  □

Our first result, without any assumption on $(M_n)$, is as follows.

THEOREM 2.1. *Let $(M_n)$ be a locally square integrable martingale. Then, for all $x, y > 0$,*

$$(2.3) \qquad \mathbb{P}(|M_n| \geq x, [M]_n + \langle M \rangle_n \leq y) \leq 2\exp\left(-\frac{x^2}{2y}\right).$$

REMARK 2.1. A similar result for continuous-time locally square integrable martingale may be found in the first part of Proposition 4.2.3 of Barlow, Jacka and Yor [5].

For self-normalized martingales, we obtain the following result.

THEOREM 2.2. *Let $(M_n)$ be a locally square integrable martingale. Then, for all $x, y > 0$, $a \geq 0$ and $b > 0$,*

$$(2.4) \quad \mathbb{P}\left(\frac{|M_n|}{a + b\langle M \rangle_n} \geq x, \langle M \rangle_n \geq [M]_n + y\right) \leq 2\exp\left(-x^2\left(ab + \frac{b^2y}{2}\right)\right).$$

*Moreover, we also have*

(2.5)
$$\mathbb{P}\left(\frac{|M_n|}{a+b\langle M\rangle_n} \geq x, [M]_n \leq y\langle M\rangle_n\right)$$
$$\leq 2\inf_{p>1}\left(\mathbb{E}\left[\exp\left(-(p-1)\frac{x^2}{(1+y)}\left(ab+\frac{b^2}{2}\langle M\rangle_n\right)\right)\right]\right)^{1/p}.$$

PROOF. The proof is given in Appendix B. □

REMARK 2.2. It is not hard to see that (2.4) and (2.5) also hold exchanging the roles of $\langle M\rangle_n$ and $[M]_n$.

**3. Random variables heavy on left or right.** This section deals with our new notion of random variables heavy on left or right. It allows us to improve Lemma 2.1.

DEFINITION 3.1. We shall say that an integrable random variable $X$ is heavy on left if $\mathbb{E}[X]=0$ and, for all $a>0$, $\mathbb{E}[T_a(X)]\leq 0$ where

$$T_a(X) = \min(|X|,a)\operatorname{sign}(X)$$

is the truncated version of $X$. Moreover, $X$ is heavy on right if $-X$ is heavy on left.

REMARK 3.1. Let $F$ be the cumulative distribution function associated with $X$. Standard calculation leads to $\mathbb{E}[T_a(X)] = -H(a)$ where $H$ is the function defined, for all $a>0$, by

$$H(a) = \int_0^a F(-x) - (1-F(x-))\,dx$$

where $F(x-)$ stands for the left limit of $F$ at point $x$. Consequently, $X$ is heavy on left if $\mathbb{E}[X]=0$ and, for all $a>0$, $H(a)\geq 0$. Moreover, $H$ is equal to zero at infinity as

$$\lim_{a\to\infty} H(a) = -\mathbb{E}[X] = 0.$$

Furthermore, one can observe that a random variable $X$ is symmetric if and only if $X$ is heavy on left and on right.

The following lemma is the keystone of our one-sided exponential inequalities.

LEMMA 3.1. *For a random variable $X$ and for all $t\in\mathbb{R}$, let*

$$L(t) = \mathbb{E}\left[\exp\left(tX - \frac{t^2}{2}X^2\right)\right].$$

(1) *If $X$ is heavy on left, then for all $t \geq 0$, $L(t) \leq 1$.*
(2) *If $X$ is heavy on right, then for all $t \leq 0$, $L(t) \leq 1$.*
(3) *If $X$ is symmetric, then for all $t \in \mathbb{R}$, $L(t) \leq 1$.*

PROOF.  The proof is given in Appendix A.  □

We shall now provide several examples of random variables heavy on left. More details concerning these examples may be found in Appendix E. We wish to point out that most of all positive random variables centered around their mean are heavy on left. As a matter of fact, let $Y$ be a positive integrable random variable with mean $m$ and denote

$$X = Y - m.$$

*Discrete random variables.*

(1) If $Y$ has a Bernoulli distribution $\mathcal{B}(p)$ with parameter $0 < p < 1$, then $X$ is heavy on left, heavy on right, or symmetric if $p < 1/2$, $p > 1/2$, or $p = 1/2$, respectively.
(2) If $Y$ has a Geometric distribution $\mathcal{G}(p)$ with parameter $0 < p < 1$, then $X$ is always heavy on left.
(3) If $Y$ has a Poisson distribution $\mathcal{P}(\lambda)$ with parameter $\lambda > 0$, then $X$ is heavy on left as soon as

$$2\exp(-\lambda) \sum_{k=0}^{[\lambda]} \frac{\lambda^k}{k!} \geq 1.$$

One can observe that this condition is always fulfilled if $\lambda$ is a positive integer; see Lemma 1 of [1].

*Continuous random variables.*

(1) If $Y$ has an exponential distribution $\mathcal{E}(\lambda)$ with parameter $\lambda > 0$, then $X$ is always heavy on left.
(2) If $Y$ has a Gamma distribution $\mathcal{G}(a, \lambda)$ with parameters $a, \lambda > 0$, then $X$ is always heavy on left.
(3) If $Y$ has a Pareto distribution with parameters $a, \lambda > 0$, that is, $Y = a\exp(Z)$ where $Z$ has an exponential distribution $\mathcal{E}(\lambda)$, then $X$ is always heavy on left.
(4) If $Y$ has a log-normal distribution with parameters $m \in \mathbb{R}$ and $\sigma^2 > 0$, that is, $Y = \exp(Z)$ where $Z$ has a Normal distribution $\mathcal{N}(m, \sigma^2)$, then $X$ is always heavy on left.

**4. One-sided exponential inequalities.** Our next results are related to martingales heavy on left in the sense of the following definition.

DEFINITION 4.1. Let $(M_n)$ be a locally square integrable martingale adapted to a filtration $\mathbb{F} = (\mathcal{F}_n)$. We shall say that $(M_n)$ is heavy on left if all its increments are conditionally heavy on left. In other words, for all $n \geq 1$ and for any $a > 0$, $\mathbb{E}[T_a(\Delta M_n)|\mathcal{F}_{n-1}] \leq 0$. Moreover, $(M_n)$ is heavy on right if $(-M_n)$ is heavy on left.

We shall recover Theorem 1.3 under the assumption that $(M_n)$ is heavy on left.

THEOREM 4.1. Let $(M_n)$ be a locally square integrable martingale heavy on left. Then, for all $x, y > 0$,

$$\mathbb{P}(M_n \geq x, [M]_n \leq y) \leq \exp\left(-\frac{x^2}{2y}\right). \tag{4.1}$$

For self-normalized martingales, our results are as follows.

THEOREM 4.2. Let $(M_n)$ be a locally square integrable martingale heavy on left. Then, for all $x > 0$, $a \geq 0$ and $b > 0$,

$$\mathbb{P}\left(\frac{M_n}{a + b[M]_n} \geq x\right) \leq \inf_{p>1}\left(\mathbb{E}\left[\exp\left(-(p-1)x^2\left(ab + \frac{b^2}{2}[M]_n\right)\right)\right]\right)^{1/p}, \tag{4.2}$$

and, for all $y > 0$,

$$\mathbb{P}\left(\frac{M_n}{a + b[M]_n} \geq x, [M]_n \geq y\right) \leq \exp\left(-x^2\left(ab + \frac{b^2 y}{2}\right)\right). \tag{4.3}$$

Moreover, we also have

$$\mathbb{P}\left(\frac{M_n}{a + b\langle M\rangle_n} \geq x, [M]_n \leq y\langle M\rangle_n\right)$$
$$\leq \inf_{p>1}\left(\mathbb{E}\left[\exp\left(-(p-1)\frac{x^2}{y}\left(ab + \frac{b^2}{2}\langle M\rangle_n\right)\right)\right]\right)^{1/p}. \tag{4.4}$$

PROOF. The proof is given in Appendix C. □

REMARK 4.1. In the particular case $p = 2$, Theorem 4.2 is due to De la Peña [8] under the conditional symmetric assumption on $(M_n)$. The only difference between (2.5) and (4.4) is that $(1 + y)$ is replaced by $y$ in the upper-bound of (4.4).

REMARK 4.2. A locally square integrable martingale $(M_n)$ is Gaussian if, for all $n \geq 1$, the distribution of its increments $\Delta M_n$ given $\mathcal{F}_{n-1}$ is $\mathcal{N}(0, \Delta \langle M \rangle_n)$. Moreover, $(M_n)$ is called sub-Gaussian if there exists some constant $\alpha > 0$ such that, for all $n \geq 1$ and $t \in \mathbb{R}$,

$$(4.5) \qquad \mathbb{E}[\exp(t\Delta M_n)|\mathcal{F}_{n-1}] \leq \exp\left(\frac{\alpha^2 t^2}{2}\Delta\langle M\rangle_n\right).$$

It is well known that if the increments of $(M_n)$ are bounded or if $(M_n)$ is Gaussian, then $(M_n)$ is sub-Gaussian. In addition, if $(M_n)$ satisfies (4.5), then inequalities (4.1), (4.2) and (4.3) hold with appropriate upper-bounds, replacing $[M]_n$ by $\langle M\rangle_n$ everywhere. For example, (4.2) can be rewritten as

$$(4.6) \qquad \begin{aligned} &\mathbb{P}\left(\frac{M_n}{a + b\langle M\rangle_n} \geq x\right) \\ &\leq \inf_{p>1}\left(\mathbb{E}\left[\exp\left(-(p-1)\frac{x^2}{\alpha^2}\left(ab + \frac{b^2}{2}\langle M\rangle_n\right)\right)\right]\right)^{1/p}. \end{aligned}$$

## 5. Applications.

5.1. *Linear regressions.* Consider the stochastic linear regression given, for all $n \geq 0$, by

$$(5.1) \qquad X_{n+1} = \theta \phi_n + \varepsilon_{n+1}$$

where $X_n$, $\phi_n$ and $\varepsilon_n$ are the observation, the regression variable and the driven noise, respectively. We assume that $(\phi_n)$ is a sequence of independent and identically distributed random variables. We also assume that $(\varepsilon_n)$ is a sequence of identically distributed random variables, with mean zero and variance $\sigma^2 > 0$. Moreover, we suppose that, for all $n \geq 0$, the random variable $\varepsilon_{n+1}$ is independent of $\mathcal{F}_n$ where $\mathcal{F}_n = \sigma(\phi_0, \varepsilon_1, \ldots, \phi_{n-1}, \varepsilon_n)$. In order to estimate the unknown parameter $\theta$, we make use of the least-squares estimator $\widehat{\theta}_n$ given, for all $n \geq 1$, by

$$(5.2) \qquad \widehat{\theta}_n = \frac{\sum_{k=1}^n \phi_{k-1} X_k}{\sum_{k=1}^n \phi_{k-1}^2}.$$

It immediately follows from (5.1) and (5.2) that

$$(5.3) \qquad \widehat{\theta}_n - \theta = \sigma^2 \frac{M_n}{\langle M\rangle_n}$$

where

$$M_n = \sum_{k=1}^n \phi_{k-1}\varepsilon_k \quad \text{and} \quad \langle M\rangle_n = \sigma^2 \sum_{k=1}^n \phi_{k-1}^2.$$

Let $H$ and $L$ be the cumulant generating functions of the sequences $(\phi_n^2)$ and $(\varepsilon_n^2)$, respectively given, for all $t \in \mathbb{R}$, by

$$H(t) = \log \mathbb{E}[\exp(t\phi_n^2)] \quad \text{and} \quad L(t) = \log \mathbb{E}[\exp(t\varepsilon_n^2)].$$

COROLLARY 5.1. *Assume that $L$ is finite on some interval $[0, c]$ with $c > 0$ and denote by $I$ its Fenchel–Legendre transform on $[0, c]$,*

$$I(x) = \sup_{0 \leq t \leq c} \{xt - L(t)\}.$$

*Then, for all $n \geq 1$, $x > 0$ and $y > 0$, we have*

(5.4)
$$\mathbb{P}(|\widehat{\theta}_n - \theta| \geq x)$$
$$\leq 2 \inf_{p>1} \exp\left(\frac{n}{p} H\left(-\frac{(p-1)x^2}{2\sigma^2(1+y)}\right)\right) + \exp\left(-nI\left(\frac{\sigma^2 y}{n}\right)\right).$$

REMARK 5.1. Corollary 5.1 is also true if $(\phi_n, \varepsilon_n)$ is a sequence of independent and identically distributed random vectors of $\mathbb{R}^2$ such that the marginal distribution of $\varepsilon_n$ is symmetric. By use of (4.4), inequality (5.4) holds replacing $(1 + y)$ by $y$ in the argument of $H$.

REMARK 5.2. As soon as the sequence $(\varepsilon_n)$ is bounded, the right-hand side of (5.4) vanishes since we may directly compare $[M_n]$ with $\langle M \rangle_n$. For example, assume that $(\varepsilon_n)$ is distributed as a centered Bernoulli $\mathcal{B}(p)$ distribution with parameter $0 < p < 1$. If $r = \max(p, q)$, we clearly have for all $n \geq 0$,

$$[M]_n \leq \frac{r^2}{pq} \langle M \rangle_n.$$

Consequently, we immediately infer from (2.5) that for all $n \geq 1$ and $x > 0$,

$$\mathbb{P}(|\widehat{\theta}_n - \theta| \geq x) \leq 2 \exp\left(\frac{n}{2} H\left(-\frac{x^2}{4r^2}\right)\right).$$

Furthermore, assume that $(\phi_n)$ is distributed as a normal $\mathcal{N}(0, \tau^2)$ distribution with variance $\tau^2 > 0$. Then, we deduce that for all $n \geq 1$ and $x > 0$,

$$\mathbb{P}(|\widehat{\theta}_n - \theta| \geq x) \leq 2 \exp\left(-\frac{n}{4} \log\left(1 + \frac{\tau^2 x^2}{2r^2}\right)\right).$$

PROOF OF COROLLARY 5.1. It follows from (2.5) that, for all $n \geq 1$, $x > 0$ and $y > 0$,

$$\mathbb{P}(|\widehat{\theta}_n - \theta| \geq x) = \mathbb{P}\left(|M_n| \geq \frac{x}{\sigma^2} \langle M \rangle_n\right) \leq P_n(x, y) + Q_n(y)$$

where $Q_n(y) = \mathbb{P}([M]_n > y\langle M\rangle_n)$ and

$$P_n(x,y) = 2\inf_{p>1}\left(\mathbb{E}\left[\exp\left(-(p-1)\frac{x^2}{2\sigma^4(1+y)}\langle M\rangle_n\right)\right]\right)^{1/p}$$

$$= 2\inf_{p>1}\exp\left(\frac{n}{p}H\left(-\frac{(p-1)x^2}{2\sigma^2(1+y)}\right)\right).$$

In addition, for all $y > 0$ and $0 \leq t \leq c$,

$$Q_n(y) \leq \mathbb{P}\left(\sum_{k=1}^n \varepsilon_k^2 > \sigma^2 y\right) \leq \exp(-\sigma^2 ty)\mathbb{E}\left[\exp\left(t\sum_{k=1}^n \varepsilon_k^2\right)\right]$$

$$\leq \exp(-\sigma^2 ty + nL(t)) \leq \exp\left(-nI\left(\frac{\sigma^2 y}{n}\right)\right),$$

which achieves the proof of Corollary 5.1. $\square$

5.2. *Autoregressive processes.* Consider the autoregressive process given, for all $n \geq 0$, by

(5.5) $$X_{n+1} = \theta X_n + \varepsilon_{n+1}$$

where $X_n$ and $\varepsilon_n$ are the observation and the driven noise, respectively. We assume that $(\varepsilon_n)$ is a sequence of independent and identically distributed random variables with standard $\mathcal{N}(0,\sigma^2)$ distribution where $\sigma^2 > 0$. The process is said to be stable if $|\theta| < 1$, unstable if $|\theta| = 1$ and explosive if $|\theta| > 1$. We can estimate the unknown parameter $\theta$ by the least-squares or the Yule–Walker estimators given, for all $n \geq 1$, by

(5.6) $$\widehat{\theta}_n = \frac{\sum_{k=1}^n X_{k-1}X_k}{\sum_{k=1}^n X_{k-1}^2} \quad \text{and} \quad \widetilde{\theta}_n = \frac{\sum_{k=1}^n X_{k-1}X_k}{\sum_{k=0}^n X_k^2}.$$

It is well known that $\widehat{\theta}_n$ and $\widetilde{\theta}_n$ both converge almost surely to $\theta$ and their fluctuations can be found in [22]. In the stable case $|\theta| < 1$, the large deviation principles were established in [6]. More precisely, set

$$a = \frac{\theta - \sqrt{\theta^2 + 8}}{4} \quad \text{and} \quad b = \frac{\theta + \sqrt{\theta^2 + 8}}{4}.$$

Assume that $X_0$ is independent of $(\varepsilon_n)$ with $\mathcal{N}(0, \sigma^2/(1-\theta^2))$ distribution. Then, $(\widehat{\theta}_n)$ and $(\widetilde{\theta}_n)$ satisfy large deviation principles with good rate functions respectively given by

$$I(x) = \begin{cases} \frac{1}{2}\log\left(\frac{1+\theta^2-2\theta x}{1-x^2}\right), & \text{if } x \in [a,b], \\ \log|\theta - 2x|, & \text{otherwise,} \end{cases}$$

$$J(x) = \begin{cases} \frac{1}{2}\log\left(\frac{1+\theta^2-2\theta x}{1-x^2}\right), & \text{if } x \in ]-1,1[, \\ +\infty, & \text{otherwise.} \end{cases}$$

It is only recently that sharp large deviation principles were established for the Yule–Walker estimator $\widetilde{\theta}_n$ in the stable, unstable and explosive cases [7]. Much work remains to be done for the least-squares estimator $\widehat{\theta}_n$. Our goal is to propose, whatever the value of $\theta$ is, a very simple exponential inequality for both $\widehat{\theta}_n$ and $\widetilde{\theta}_n$. For the sake of simplicity, we assume that $X_0$ is independent of $(\varepsilon_n)$ with $\mathcal{N}(0, \tau^2)$ distribution where $\tau^2 \geq \sigma^2$.

COROLLARY 5.2.　*For all $n \geq 1$ and $x > 0$, we have*

$$\text{(5.7)} \qquad \mathbb{P}(|\widehat{\theta}_n - \theta| \geq x) \leq 2 \exp\left(-\frac{nx^2}{2(1+y_x)}\right)$$

*where $y_x$ is the unique positive solution of the equation $h(y_x) = x^2$ and $h$ is the function $h(x) = (1+x)\log(1+x) - x$. Moreover, for all $n \geq 1$ and $x > 0$, we also have*

$$\text{(5.8)} \qquad \mathbb{P}(|\widetilde{\theta}_n - \theta| \geq x + |\theta|) \leq 2 \exp\left(-\frac{nx^2}{2(1+y_x)}\right).$$

PROOF.　The proof is given in Appendix D.　□

REMARK 5.3.　Inequality (5.7) can be very simple if $x$ is small enough. As a matter of fact, one can easily see that for all $0 < x < 1$, $h(x) < x^2/4$. Consequently, it immediately follows from (5.7) that, for all $0 < x < 1/2$,

$$\mathbb{P}(|\widehat{\theta}_n - \theta| \geq x) \leq 2 \exp\left(-\frac{nx^2}{2(1+2x)}\right).$$

Moreover, if $\theta > 0$, we can deduce from (5.6) that, for all $x > 0$,

$$\mathbb{P}(\widetilde{\theta}_n - \theta \geq x) \leq \exp\left(-\frac{nx^2}{2(1+y_x)}\right).$$

5.3. *Branching processes.*　Consider the Galton–Watson process starting from $X_0 = 1$ and given, for all $n \geq 1$, by

$$\text{(5.9)} \qquad X_n = \sum_{k=1}^{X_{n-1}} Y_{n,k}$$

where $(Y_{n,k})$ is a sequence of independent and identically distributed, nonnegative integer-valued random variables. The distribution of $(Y_{n,k})$, with finite mean $m$ and variance $\sigma^2$, is commonly called the offspring or reproduction distribution. Hereafter, we shall assume that $m > 1$. In order to estimate the offspring mean $m$, we can make use of the Lotka–Nagaev or the Harris estimators given, for all $n \geq 1$, by

$$\text{(5.10)} \qquad \widetilde{m}_n = \frac{X_n}{X_{n-1}} \quad \text{and} \quad \widehat{m}_n = \frac{\sum_{k=1}^n X_k}{\sum_{k=1}^n X_{k-1}}.$$

Without loss of generality, we can suppose that the set of extinction of the process $(X_n)$ is negligible. Consequently, the Lotka–Nagaev estimator $\widetilde{m}_n$ is always well defined. It is well known that $\widetilde{m}_n$ and $\widehat{m}_n$ both converge almost surely to $m$ and their fluctuations are given in [3, 14, 15]. Moreover, the large deviation properties associated with $(\widetilde{m}_n)$ may be found in [2, 19, 20]. Our goal is now to establish, as in the previous sections, exponential inequalities for both $\widetilde{m}_n$ and $\widehat{m}_n$. Denote by $L$ the cumulant generating function associated with the centered offspring distribution given, for all $t \in \mathbb{R}$, by $L(t) = \log \mathbb{E}[\exp(t(Y_{n,k} - m))]$.

COROLLARY 5.3. *Assume that $L$ is finite on some interval $[-c, c]$ with $c > 0$ and let $I$ be its Fenchel–Legendre transform,*

$$I(x) = \sup_{-c \leq t \leq c} \{xt - L(t)\}.$$

*Then, for all $n \geq 1$ and $x > 0$,*

(5.11) $$\mathbb{P}(|\widetilde{m}_n - m| \geq x) \leq 2\mathbb{E}[\exp(-J(x)X_{n-1})]$$

*where $J(x) = \min(I(x), I(-x))$. Moreover, we also have*

(5.12) $$\mathbb{P}(|\widetilde{m}_n - m| \geq x) \leq 2\inf_{p>1}(\mathbb{E}[\exp(-(p-1)J(x)X_{n-1})])^{1/p}.$$

*In addition, if $S_n = \sum_{k=0}^n X_k$, we have for all $n \geq 1$ and $x > 0$,*

(5.13) $$\mathbb{P}(|\widehat{m}_n - m| \geq x) \leq 2\inf_{p>1}(\mathbb{E}[\exp(-(p-1)J(x)S_{n-1})])^{1/p}.$$

PROOF. The proof is given in Appendix E. □

REMARK 5.4. On the one hand, inequality (5.12) obviously holds for the Harris estimator $\widehat{m}_n$ since we always have $S_n \geq X_n$. On the other hand, in order to specify the right-hand side of (5.11), (5.12) or (5.13), it is necessary to find an upper-bound or to provide an explicit expression of the moment generating function of $X_n$. One can easily carry out this calculation when the offspring distribution is the geometric $\mathcal{G}(p)$ distribution with parameter $0 < p < 1$. As a matter of fact, in that particular case, the offspring mean $m = 1/p$ and it follows from formula (7.3) of [15] that for all $0 < s < 1$,

$$\mathbb{E}[s^{X_n}] \leq \frac{p^n s}{1 - s}.$$

Consequently, for all $n \geq 1$ and $x > 0$, we obtain the simple inequality

$$\mathbb{P}(|\widetilde{m}_n - m| \geq x) \leq \frac{2p^n \exp(-J(x))}{p(1 - \exp(-J(x)))}.$$

If the offspring distribution is not geometric, one can precisely estimate the moment generating function of $X_n$ using Theorem 1, page 80 of [3] which gives a good approximation of the distribution of $X_n$ based on the limiting distribution

$$W = \lim_{n \to \infty} \frac{X_n}{m^n} \quad \text{a.s.}$$

## APPENDIX A

This appendix is devoted to the proofs of Lemma 2.1 and Lemma 3.1. Lemma 2.1 immediately follows from Jensen's inequality. As a matter of fact, (2.1) implies that for all $t \in \mathbb{R}$,

$$L(t) \geq \exp\left(\mathbb{E}\left[tX - \frac{t^2}{2}X^2\right]\right) \geq \exp\left(-\frac{t^2}{2}\sigma^2\right).$$

Consequently, we obtain that for all $t \in \mathbb{R}$,

$$(A.1) \qquad L(t) \geq 1 - \frac{t^2}{2}\sigma^2.$$

Furthermore, for all $t \in \mathbb{R}$,

$$(A.2) \qquad L(t) + L(-t) = 2\mathbb{E}\left[\exp\left(-\frac{t^2}{2}X^2\right)\cosh(tX)\right] \leq 2$$

by the well-known inequality $\cosh(x) \leq \exp(x^2/2)$. Hence, we obtain from (A.1) together with (A.2) that for all $t \in \mathbb{R}$,

$$L(t) \leq 2 - L(-t) \leq 1 + \frac{t^2}{2}\sigma^2.$$

Lemma 3.1 is much more difficult to prove. Let $f$ be the function defined, for all $x \in \mathbb{R}$, by

$$f(x) = \exp\left(x - \frac{x^2}{2}\right).$$

We clearly have $f'(x) = (1-x)f(x)$ and $f'(-x) = (1+x)f(-x)$. We shall also make use of the functions $a$ and $b$ defined, for all $x \in \mathbb{R}$, by $a(x) = f'(-x)$ and $b(x) = f'(-x) - f'(x)$. One can realize that, for all $x > 0$, $0 < a(x) < 1$, $0 < b(x) < 2$ and $a'(x) < 0$ as $a'(x) = -(2x + x^2)f(-x)$. After those simple preliminaries, we are in position to prove Lemma 3.1. For all $t \in \mathbb{R}$,

$$L(t) = \mathbb{E}\left[\exp\left(tX - \frac{t^2X^2}{2}\right)\right] = \int_{\mathbb{R}} f(tx)\,dF(x)$$

where $F$ is the distribution function associated with $X$. Integrating by parts, we have for all $t \in \mathbb{R}$,

$$L(t) = -t \int_{\mathbb{R}} f'(tx) F(x) \, dx$$

(A.3)
$$= -t \int_{-\infty}^{0} f'(tx) F(x) \, dx - t \int_{0}^{+\infty} f'(tx) F(x) \, dx$$

$$= -t \int_{0}^{+\infty} f'(-tx) F(-x) \, dx - t \int_{0}^{+\infty} f'(tx) F(x) \, dx.$$

Consequently, as

$$-t \int_{0}^{\infty} f'(tx) \, dx = \left[ -\exp\left(tx - \frac{t^2 x^2}{2}\right) \right]_{0}^{+\infty} = 1,$$

we obtain from (A.3) that, for all $t \in \mathbb{R}$, $L(t) = 1 - tI(t)$ where

(A.4)
$$I(t) = \int_{0}^{+\infty} f'(-tx) F(-x) \, dx - \int_{0}^{+\infty} f'(tx)(1 - F(x)) \, dx$$
$$= A(t) + B(t)$$

with

$$A(t) = \int_{0}^{+\infty} a(tx)(F(-x) - (1 - F(x))) \, dx,$$

$$B(t) = \int_{0}^{+\infty} b(tx)(1 - F(x)) \, dx.$$

First of all, assume that $X$ is heavy on left. Our goal is to show that, for all $t > 0$, the integral $I(t)$ is nonnegative. We obviously have, for all $t > 0$, $B(t) \geq 0$ as $b(tx) > 0$. In addition, for any $a > 0$, let

$$H(a) = \int_{0}^{a} F(-x) - (1 - F(x-)) \, dx.$$

Since $H'(a) = F(-a) - (1 - F(a))$ almost everywhere, integrating once again by parts, we find that

(A.5)
$$A(t) = \int_{0}^{+\infty} a(tx) H'(x) \, dx = [a(tx) H(x)]_{0}^{+\infty} - \int_{0}^{+\infty} ta'(tx) H(x) \, dx$$
$$= -t \int_{0}^{+\infty} a'(tx) H(x) \, dx$$

as $H(0) = 0$ and $H$ vanishes at infinity. Hereafter, as $X$ is heavy on left, $H(a) \geq 0$ for all $a \geq 0$. Moreover, we recall that, for all $x > 0$, $a'(x) < 0$. Hence, we immediately deduce from (A.5) that, for all $t > 0$, $A(t) \geq 0$. Consequently, relation (A.4) leads to $I(t) \geq 0$ and $L(t) \leq 1$ for all $t > 0$, which completes the proof of part (1) of Lemma 3.1. Next, if $X$ is heavy on right, $-X$ is heavy on left. Hence, we immediately infer from (2.1) and part (1) of Lemma 3.1 that $L(t) \leq 1$ for all $t < 0$. Finally, part (3) of Lemma 3.1 follows from the conjunction of parts (1) and (2). Another straightforward way to prove part (3) is as follows. If $X$ is symmetric, we have for all $t \in \mathbb{R}$,

$$L(t) = \int_{\mathbb{R}} f(tx) \, dF(x) = \int_0^{+\infty} (f(tx) + f(-tx)) \, dF(x)$$
$$= 2 \int_0^{+\infty} \exp(-t^2 x^2/2) \cosh(tx) \, dF(x) \leq 1$$

by the well-known inequality $\cosh(x) \leq \exp(x^2/2)$.

## APPENDIX B

In order to prove Theorems 2.1 and 2.2, we shall often make use of the following lemma.

LEMMA B.1. *Let $(M_n)$ be a locally square integrable martingale. For all $t \in \mathbb{R}$ and $n \geq 0$, denote*

$$V_n(t) = \exp\left(tM_n - \frac{t^2}{2}([M]_n + \langle M \rangle_n)\right).$$

*Then, for all $t \in \mathbb{R}$, $(V_n(t))$ is a positive supermartingale with $\mathbb{E}[V_n(t)] \leq 1$.*

PROOF. For all $t \in \mathbb{R}$ and $n \geq 1$, we have

$$V_n(t) = V_{n-1}(t) \exp\left(t\Delta M_n - \frac{t^2}{2}(\Delta[M]_n + \Delta\langle M \rangle_n)\right)$$

where $\Delta M_n = M_n - M_{n-1}$, $\Delta[M]_n = \Delta M_n^2$ and $\Delta \langle M \rangle_n = \mathbb{E}[\Delta M_n^2 | \mathcal{F}_{n-1}]$. Hence, we deduce from Lemma 2.1 that for all $t \in \mathbb{R}$,

$$\mathbb{E}[V_n(t)|\mathcal{F}_{n-1}] \leq V_{n-1}(t) \exp\left(-\frac{t^2}{2}\Delta\langle M \rangle_n\right)\left(1 + \frac{t^2}{2}\Delta\langle M \rangle_n\right)$$
$$\leq V_{n-1}(t).$$

Consequently, for all $t \in \mathbb{R}$, $(V_n(t))$ is a positive supermartingale such that, for all $n \geq 1$, $\mathbb{E}[V_n(t)] \leq \mathbb{E}[V_{n-1}(t)]$ which implies that $\mathbb{E}[V_n(t)] \leq \mathbb{E}[V_0(t)] = 1$. □

We are now in position to prove Theorems 2.1 and 2.2 inspired by the original work of De la Peña [8]. First of all, denote

$$Z_n = [M]_n + \langle M \rangle_n.$$

For all $x, y > 0$, let
$$A_n = \{|M_n| \geq x, Z_n \leq y\}.$$
We have the decomposition $A_n = A_n^+ \cup A_n^-$ where $A_n^+ = \{M_n \geq x, Z_n \leq y\}$ and $A_n^- = \{M_n \leq -x, Z_n \leq y\}$. By Markov's inequality, we have for all $t > 0$,

$$\mathbb{P}(A_n^+) \leq \mathbb{E}\left[\exp\left(\frac{t}{2}M_n - \frac{tx}{2}\right)\mathbb{1}_{A_n^+}\right]$$
$$\leq \mathbb{E}\left[\exp\left(\frac{t}{2}M_n - \frac{t^2}{4}Z_n\right)\exp\left(\frac{t^2}{4}Z_n - \frac{tx}{2}\right)\mathbb{1}_{A_n^+}\right]$$
$$\leq \exp\left(\frac{t^2 y}{4} - \frac{tx}{2}\right)\sqrt{\mathbb{E}[V_n(t)]\mathbb{P}(A_n^+)}.$$

Hence, we deduce from Lemma B.1 that for all $t > 0$,

(B.1) $$\mathbb{P}(A_n^+) \leq \exp\left(\frac{t^2 y}{4} - \frac{tx}{2}\right)\sqrt{\mathbb{P}(A_n^+)}.$$

Dividing both sides of (B.1) by $\sqrt{\mathbb{P}(A_n^+)}$ and choosing the value $t = x/y$, we find that
$$\mathbb{P}(A_n^+) \leq \exp\left(-\frac{x^2}{2y}\right).$$
We also find the same upper-bound for $\mathbb{P}(A_n^-)$ which immediately leads to (2.3).

We next proceed to the proof of Theorem 2.2 in the special case $a = 0$ and $b = 1$ inasmuch as the proof for the general case follows exactly the same lines. For all $x, y > 0$, let
$$B_n = \{|M_n| \geq x\langle M\rangle_n, \langle M\rangle_n - [M]_n \geq y\} = B_n^+ \cup B_n^-$$
where
$$B_n^+ = \{M_n \geq x\langle M\rangle_n, \langle M\rangle_n - [M]_n \geq y\},$$
$$B_n^- = \{M_n \leq -x\langle M\rangle_n, \langle M\rangle_n - [M]_n \geq y\}.$$
By Cauchy–Schwarz's inequality, we have for all $t > 0$,

(B.2) $$\mathbb{P}(B_n^+) \leq \mathbb{E}\left[\exp\left(\frac{t}{2}M_n - \frac{tx}{2}\langle M\rangle_n\right)\mathbb{1}_{B_n^+}\right]$$
$$\leq \mathbb{E}\left[\exp\left(\frac{t}{2}M_n - \frac{t^2}{4}Z_n\right)\exp\left(\frac{t}{4}(t-2x)\langle M\rangle_n + \frac{t^2}{4}[M]_n\right)\mathbb{1}_{B_n^+}\right].$$

Consequently, we obtain from (B.2) with the particular choice $t = x$ that

(B.3) $$\mathbb{P}(B_n^+) \leq \exp\left(-\frac{x^2 y}{4}\right)\sqrt{\mathbb{P}(B_n^+)}.$$

Therefore, if we divide both sides of (B.3) by $\sqrt{\mathbb{P}(B_n^+)}$, we find that

$$\mathbb{P}(B_n^+) \leq \exp\left(-\frac{x^2 y}{2}\right).$$

The same upper-bound holds for $\mathbb{P}(B_n^-)$ which clearly implies (2.4). Furthermore, for all $x, y > 0$, let

$$C_n = \{|M_n| \geq x\langle M\rangle_n, [M]_n \leq y\langle M\rangle_n\} = C_n^+ \cup C_n^-$$

where

$$C_n^+ = \{M_n \geq x\langle M\rangle_n, [M]_n \leq y\langle M\rangle_n\},$$
$$C_n^- = \{M_n \leq -x\langle M\rangle_n, [M]_n \leq y\langle M\rangle_n\}.$$

By Hölder's inequality, we have for all $t > 0$ and $q > 1$,

$$\mathbb{P}(C_n^+) \leq \mathbb{E}\left[\exp\left(\frac{t}{q}M_n - \frac{tx}{q}\langle M\rangle_n\right)\mathbb{1}_{C_n^+}\right]$$

(B.4)
$$\leq \mathbb{E}\left[\exp\left(\frac{t}{q}M_n - \frac{t^2}{2q}Z_n\right)\exp\left(\frac{t}{2q}(t - 2x + ty)\langle M\rangle_n\right)\mathbb{1}_{C_n^+}\right]$$

$$\leq \left(\mathbb{E}\left[\exp\left(\frac{tp}{2q}(t - 2x + ty)\langle M\rangle_n\right)\right]\right)^{1/p}.$$

Consequently, as $p/q = p - 1$, we can deduce from (B.4) and the particular choice $t = x/(1 + y)$ that

$$\mathbb{P}(C_n^+) \leq \inf_{p>1}\left(\mathbb{E}\left[\exp\left(-(p-1)\frac{x^2}{2(1+y)}\langle M\rangle_n\right)\right]\right)^{1/p}.$$

We also find the same upper-bound for $\mathbb{P}(C_n^-)$ which completes the proof of Theorem 2.2.

## APPENDIX C

The proofs of Theorems 4.1 and 4.2 are based on the following lemma.

LEMMA C.1. *Let $(M_n)$ be a locally square integrable martingale. For all $t \in \mathbb{R}$ and $n \geq 0$, denote*

$$W_n(t) = \exp\left(tM_n - \frac{t^2}{2}[M]_n\right).$$

(1) *If $(M_n)$ is heavy on left, then for all $t \geq 0$, $(W_n(t))$ is a supermartingale with $\mathbb{E}[W_n(t)] \leq 1$.*
(2) *If $(M_n)$ is heavy on right, then for all $t \leq 0$, $(W_n(t))$ is a supermartingale with $\mathbb{E}[W_n(t)] \leq 1$.*

(3) If $(M_n)$ is conditionally symmetric, then for all $t \in \mathbb{R}$, $(W_n(t))$ is a supermartingale with $\mathbb{E}[W_n(t)] \leq 1$.

PROOF. Lemma C.1 part (3) is due to De la Peña [8], Lemma 6.1. Our approach is totally different as it mainly relies on Lemma 3.1. Assume that $(M_n)$ is heavy on left. For all $t \in \mathbb{R}$ and $n \geq 1$, we have

$$W_n(t) = W_{n-1}(t) \exp\left(t \Delta M_n - \frac{t^2}{2} \Delta[M]_n\right)$$

where $\Delta[M]_n = \Delta M_n^2$. We infer from Lemma 3.1 part (1) that for all $n \geq 1$ and for all $t \geq 0$,

$$\mathbb{E}\left[\exp\left(t \Delta M_n - \frac{t^2}{2} \Delta M_n^2\right) \Big| \mathcal{F}_{n-1}\right] \leq 1.$$

Consequently, for all $t \geq 0$, $(W_n(t))$ is a positive supermartingale such that, for all $n \geq 1$, $\mathbb{E}[W_n(t)] \leq \mathbb{E}[W_{n-1}(t)]$ which leads to $\mathbb{E}[W_n(t)] \leq \mathbb{E}[W_0(t)] = 1$. The rest of the proof is also a straightforward application of Lemma 3.1. □

By use of Lemma C.1, the proof of Theorem 4.1 is quite analogous to that of Theorem 2.1 and therefore is left to the reader. We shall proceed to the proof of Theorem 4.2 in the special case $a = 0$ and $b = 1$. For all $x > 0$, let $A_n = \{M_n \geq x[M]_n\}$. By Hölder's inequality, we have for all $t > 0$ and $q > 1$,

$$\mathbb{P}(A_n) \leq \mathbb{E}\left[\exp\left(\frac{t}{q} M_n - \frac{tx}{q}[M]_n\right) \mathbb{1}_{A_n}\right]$$

(C.1)
$$\leq \mathbb{E}\left[\exp\left(\frac{t}{q} M_n - \frac{t^2}{2q}[M]_n\right) \exp\left(\frac{t}{2q}(t - 2x)[M]_n\right) \mathbb{1}_{A_n}\right]$$

$$\leq \left(\mathbb{E}\left[\exp\left(\frac{tp}{2q}(t - 2x)[M]_n\right)\right]\right)^{1/p} (\mathbb{E}[W_n(t)])^{1/q}.$$

Since $(M_n)$ is heavy on left, it follows from Lemma C.1 that for all $t \geq 0$, $\mathbb{E}[W_n(t)] \leq 1$. Consequently, as $p/q = p - 1$, we can deduce from (C.1) and the particular choice $t = x$ that

$$\mathbb{P}(A_n) \leq \inf_{p > 1} \left(\mathbb{E}\left[\exp\left(-(p-1)\frac{x^2}{2}[M]_n\right)\right]\right)^{1/p}.$$

Furthermore, for all $x, y > 0$, let $B_n = \{M_n \geq x[M]_n, [M]_n \geq y\}$. As before, we find that for all $0 < t < 2x$,

$$\mathbb{P}(B_n) \leq \mathbb{E}\left[\exp\left(\frac{t}{2} M_n - \frac{t^2}{4}[M]_n\right) \exp\left(\frac{t}{4}(t - 2x)[M]_n\right) \mathbb{1}_{B_n}\right]$$

$$\leq \exp\left(\frac{ty}{4}(t - 2x)\right) \mathbb{E}\left[\exp\left(\frac{t}{2} M_n - \frac{t^2}{4}[M]_n\right) \mathbb{1}_{B_n}\right]$$

$$\leq \exp\left(\frac{ty}{4}(t-2x)\right)\sqrt{\mathbb{P}(B_n)}$$

$$\leq \exp\left(-\frac{x^2y}{2}\right),$$

choosing the value $t = x$. Finally, the last inequality of Theorem 4.2 is left to the reader as its proof follows exactly the same arguments as (4.3).

## APPENDIX D

We shall now focus our attention on the proof of Corollary 5.2. It immediately follows from (5.5) together with (5.6) that for all $n \geq 1$,

(D.1) $$\widehat{\theta}_n - \theta = \sigma^2 \frac{M_n}{\langle M \rangle_n}$$

where

$$M_n = \sum_{k=1}^n X_{k-1}\varepsilon_k \quad \text{and} \quad \langle M \rangle_n = \sigma^2 \sum_{k=1}^n X_{k-1}^2.$$

The driven noise $(\varepsilon_n)$ is a sequence of independent and identically distributed random variables with $\mathcal{N}(0, \sigma^2)$ distribution. Consequently, for all $n \geq 1$, the distribution of the increments $\Delta M_n = X_{n-1}\varepsilon_n$ given $\mathcal{F}_{n-1}$ is $\mathcal{N}(0, \sigma^2 X_{n-1}^2)$ which implies that $(M_n)$ is a Gaussian martingale. Therefore, we infer from inequality (4.6) that for all $n \geq 1$ and $x > 0$,

(D.2) $$\mathbb{P}(|\widehat{\theta}_n - \theta| \geq x) = \mathbb{P}\left(|M_n| \geq \frac{x}{\sigma^2}\langle M \rangle_n\right) = 2\mathbb{P}\left(M_n \geq \frac{x}{\sigma^2}\langle M \rangle_n\right)$$

$$\leq 2 \inf_{p>1}\left(\mathbb{E}\left[\exp\left(-(p-1)\frac{x^2}{2\sigma^4}\langle M \rangle_n\right)\right]\right)^{1/p}.$$

Similar result may be found in [17, 23, 24]. We are now halfway to our goal and it remains to find a suitable upper-bound for the right-hand side of (D.2). For all $t \in \mathbb{R}$ such that $1 - 2\sigma^2 t > 0$, if $\alpha = 1/\sqrt{1 - 2\sigma^2 t}$, we deduce from (5.5) that, for all $n \geq 1$,

$$\mathbb{E}[\exp(tX_n^2)|\mathcal{F}_{n-1}] = \exp(t\theta^2 X_{n-1}^2)\mathbb{E}[\exp(2\theta t X_{n-1}\varepsilon_n + t\varepsilon_n^2)|\mathcal{F}_{n-1}]$$

$$= \frac{\exp(t\theta^2 X_{n-1}^2)}{\sigma\sqrt{2\pi}} \int_\mathbb{R} \exp\left(-\frac{x^2}{2\alpha^2\sigma^2}\right) \exp(2\theta t X_{n-1} x)\, dx.$$

Hence, if $\beta = 2t\alpha\sigma\theta X_{n-1}$, we find via the change of variables $y = x/\alpha\sigma$ that

$$\mathbb{E}[\exp(tX_n^2)|\mathcal{F}_{n-1}] = \frac{\alpha \exp(t\theta^2 X_{n-1}^2)}{\sqrt{2\pi}} \int_\mathbb{R} \exp\left(-\frac{y^2}{2} + \beta y\right) dy$$

$$= \alpha \exp\left(t\theta^2 X_{n-1}^2 + \frac{\beta^2}{2}\right) = \alpha \exp(t\alpha^2\theta^2 X_{n-1}^2),$$

which implies that, for all $t < 0$ and $n \geq 1$,

(D.3) $$\mathbb{E}[\exp(tX_n^2)|\mathcal{F}_{n-1}] \leq \alpha.$$

Furthermore, as $X_0$ is $\mathcal{N}(0, \tau^2)$ distributed with $\tau^2 \geq \sigma^2$, $\mathbb{E}[\exp(tX_0^2)] \leq \alpha$. It immediately follows from (D.3) together with the tower property of the conditional expectation that for all $t < 0$ and $n \geq 0$,

(D.4) $$\mathbb{E}[\exp(t\langle M\rangle_n)] \leq (1 - 2\sigma^4 t)^{-n/2}.$$

Consequently, we deduce from the conjunction of (D.2) and (D.4) with the value $t = -(p-1)x^2/2\sigma^4$ and the change of variables $y = (p-1)x^2$ that for all $x > 0$ and $n \geq 1$,

$$\mathbb{P}(|\widehat{\theta}_n - \theta| \geq x) \leq 2 \inf_{y > 0} \exp\left(-\frac{nx^2}{2}\ell(y)\right)$$

where the function $\ell$ is given by

$$\ell(y) = \frac{\log(1+y)}{x^2 + y}.$$

We clearly have

$$\ell'(y) = \frac{x^2 - h(y)}{(1+y)(x^2+y)^2}$$

where $h(y) = (1+y)\log(1+y) - y$. One can observe that the function $h$ is the Cramer transform of the centered Poisson distribution with parameter 1. Let $y_x$ be the unique positive solution of the equation $h(y_x) = x^2$. The value $y_x$ maximizes the function $\ell$ and this natural choice clearly leads to (5.7). Finally, it follows from (5.6) and (5.7) that for all $x > 0$ and $n \geq 1$,

(D.5) $$\mathbb{P}(|\widetilde{\theta}_n - \theta + \theta f_n| \geq x) \leq 2 \exp\left(-\frac{nx^2}{2(1+y_x)}\right)$$

where the random variable $0 \leq f_n \leq 1$. Hence, (D.5) implies (5.8) which completes the proof of Corollary 5.2.

## APPENDIX E

We shall now proceed to the proof of Corollary 5.3. We only focus our attention on the Harris estimator inasmuch as the proof for the Lotka–Nagaev estimator follows essentially the same lines. First of all, relation (5.9) can be rewritten as

(E.1) $$X_n = mX_{n-1} + \xi_n$$

where $\xi_n = X_n - \mathbb{E}[X_n|\mathcal{F}_{n-1}]$. Consequently, we obtain from (5.10) together with (E.1) that for all $n \geq 1$,

$$\widehat{m}_n - m = \frac{M_n}{S_{n-1}} \qquad \text{where } M_n = \sum_{k=1}^n \xi_k.$$

Moreover, for all $n \geq 1$ and $0 \leq t \leq c$, $\mathbb{E}[\exp(t\xi_n)|\mathcal{F}_{n-1}] = \exp(X_{n-1}L(t))$ which implies that

(E.2) $$\mathbb{E}[\exp(tM_n - L(t)S_{n-1})] = 1.$$

We are in position to prove (5.13). For all $x > 0$, let $D_n = \{|\widehat{m}_n - m| \geq x\}$. We have the decomposition $D_n = D_n^+ \cup D_n^-$ where $D_n^+ = \{\widehat{m}_n - m \geq x\}$ and $D_n^- = \{\widehat{m}_n - m \leq -x\}$. By Hölder's inequality together with (E.2), we have for all $0 \leq t \leq c$ and $q > 1$,

$$\mathbb{P}(D_n^+) \leq \mathbb{E}\left[\exp\left(\frac{t}{q}M_n - \frac{tx}{q}S_{n-1}\right)\mathbb{1}_{D_n^+}\right]$$

(E.3) $$\leq \mathbb{E}\left[\exp\left(\frac{t}{q}M_n - \frac{L(t)}{q}S_{n-1}\right)\exp\left(\frac{1}{q}(L(t) - tx)S_{n-1}\right)\mathbb{1}_{D_n^+}\right]$$

$$\leq \left(\mathbb{E}\left[\exp\left(\frac{p}{q}(L(t) - tx)S_{n-1}\right)\right]\right)^{1/p}.$$

Taking the infimum over the interval $[0, c]$, we infer from (E.3) that

(E.4) $$\mathbb{P}(D_n^+) \leq (\mathbb{E}[\exp(-(p-1)I(x)S_{n-1})])^{1/p}.$$

Along the same lines, we also find that

(E.5) $$\mathbb{P}(D_n^-) \leq (\mathbb{E}[\exp(-(p-1)I(-x)S_{n-1})])^{1/p}.$$

Finally, (5.13) immediately follows from (E.4) and (E.5).

## APPENDIX F

This appendix is devoted to some justifications about the examples of random variables heavy on left or right. Consider an integrable random variable $X$ with zero mean and denote by $F$ its cumulative distribution function. Let $H$ be the function defined, for all $a > 0$, by

$$H(a) = \int_0^a F(-x) - (1 - F(x-))\,dx.$$

We already saw that $X$ is heavy on left if, for all $a > 0$, $H(a) \geq 0$ while $X$ is heavy on right if, for all $a > 0$, $H(a) \leq 0$. Let $Y$ be a positive integrable random variable with mean $m$ and denote

$$X = Y - m.$$

*Discrete random variables.* Assume that $Y$ is a discrete random variable taking its values in $\mathbb{N}$. For all $n \geq 0$, let

$$s_n = \sum_{k=0}^{n} \mathbb{P}(Y = k).$$

After some straightforward calculations, we obtain that, for all $a > 0$,

$$H(a) = -a + \sum_{k=[m-a]}^{[m+a]} s_k - s_{[m+a]} + \{m+a\}s_{[m+a]} - \{m-a\}s_{[m-a]}$$

where, for all $x \in \mathbb{R}$, $[x]$ stands for the integer part of $x$ and its fractional part $\{x\}$ is given by $\{x\} = x - [x]$ and, of course, $s_n = 0$ for all $n < 0$.

*Continuous random variables.* Assume that $Y$ is a real random variable absolutely continuous with respect to the Lebesgue measure. Denote by $g$ its probability density function. It is not hard to see that, for all $a > 0$,

$$H(a) = -a + 2a \int_0^{a_m} g(x)\,dx + \int_{a_m}^{m+a} (m+a-x)g(x)\,dx$$

where $a_m = \inf\{m-a, 0\}$. Consequently, in order to check that $X$ is heavy on left, it is only necessary to show that, for all $a > 0$, $H(a) \geq 0$.

**Acknowledgments.** The authors would like to thank V. De la Peña, T. L. Lai and M. Ledoux for fruitful discussions.

## REFERENCES

[1] ADELL, J. A. and JODRÁ, P. (2005). The median of the Poisson distribution. *Metrika* **61** 337–346. MR2230380
[2] ATHREYA, K. B. (1994). Large deviation rates for branching processes I. Single type case. *Ann. Appl. Probab.* **4** 779–790. MR1284985
[3] ATHREYA, K. B. and NEY, P. E. (1972). *Branching Processes.* Springer, Berlin. MR0373040
[4] AZUMA, K. (1967). Weighted sums of certain dependent random variables. *Tôkuku Math. J.* **19** 357–367. MR0221571
[5] BARLOW, M. T., JACKA, S. D. and YOR, M. (1986). Inequalities for a pair of processes stopped at a random time. *Proc. London Math. Soc.* **52** 142–172. MR0812449
[6] BERCU, B., GAMBOA, F. and ROUAULT, A. (1997). Large deviations for quadratic forms of stationary Gaussian processes. *Stochastic Process. Appl.* **71** 75–90. MR1480640
[7] BERCU, B. (2001). On large deviations in the Gaussian autoregressive process: Stable, unstable and explosive cases. *Bernoulli* **7** 299–316. MR1828507
[8] DE LA PEÑA, V. H. (1999). A general class of exponential inequalities for martingales and ratios. *Ann. Probab.* **27** 537–564. MR1681153
[9] DE LA PEÑA, V. H., KLASS, M. J. and LAI, T. L. (2004). Self-normalized processes: Exponential inequalities, moments bounds and iterated logarithm law. *Ann. Probab.* **32** 1902–1933. MR2073181

[10] DE LA PEÑA, V. H., KLASS, M. J. and LAI, T. L. (2007). Pseudo-maximization and self-normalized processes. *Probab. Surveys* **4** 172–192. MR2368950

[11] DZHAPARIDZE, K. and VAN ZANTEN, J. H. (2001). On Bernstein-type inequalities for martingales. *Stochastic Process. Appl.* **93** 109–117. MR1819486

[12] VAN DE GEER, S. (1995). Exponential inequalities for martingales, with application to maximum likelihood estimation for counting processes. *Ann. Statist.* **23** 1779–1801. MR1370307

[13] FREEDMAN, D. A. (1975). On tail probabilities for martingales. *Ann. Probab.* **3** 100–118. MR0380971

[14] GUTTORP, P. (1991). *Statistical Inference for Branching Processes*. Wiley, New York. MR1254434

[15] HARRIS, T. E. (1963). *The Theory of Branching Processes*. Springer, Berlin. MR0163361

[16] HOEFFDING, W. J. (1963). Probability inequalities sums of bounded random variables. *J. Amer. Statist. Assoc.* **58** 713–721. MR0144363

[17] LIPTSER, R. and SPOKOINY, V. (2000). Deviation probability bound for martingales with applications to statistical estimation. *Statist. Probab. Lett.* **46** 347–357. MR1743992

[18] MCDIARMID, C. (1998). Concentration. In *Probabilistic Methods for Algorithmic Discrete Mathematics* (M. Habib, C. McDiarmid, T. Ramirez-Alfonsin and B. Reed, eds.) 195–248. Springer, Berlin. MR1678578

[19] NEY, P. E. and VIDYASHANKAR, A. N. (2003). Harmonic moments and large deviation rates for supercritical branching processes. *Ann. Appl. Probab.* **13** 475–489. MR1970272

[20] NEY, P. E. and VIDYASHANKAR, A. N. (2004). Local limit theory and large deviations for supercritical branching processes. *Ann. Appl. Probab.* **14** 1135–1166. MR2071418

[21] PINELIS, I. (1994). Optimum bounds for the distributions of martingales in Banach spaces. *Ann. Probab.* **22** 1679–1706. MR1331198

[22] WHITE, J. S. (1958). The limit distribution of the serial correlation in the explosive case. *Ann. Math. Statist.* **29** 1188–1197. MR0100952

[23] WORMS, J. (1999). Moderate deviations for stable Markov chains and regression models. *Electron. J. Probab.* **4** 1–28. MR1684149

[24] WORMS, J. (2001). Large and moderate deviations upper bounds for the Gaussian autoregressive process. *Statist. Probab. Lett.* **51** 235–243. MR1822730

INSTITUT DE MATHÉMATIQUES DE BORDEAUX
UNIVERSITÉ BORDEAUX 1
UMR 5251
351 COURS DE LA LIBÉRATION
33405 TALENCE CEDEX
FRANCE
E-MAIL: Bernard.Bercu@math.u-bordeaux1.fr

DÉPARTEMENT DE MATHÉMATIQUES
FACULTÉ DES SCIENCES DE BIZERTE
7021 ZARZOUNA
TUNISIE
E-MAIL: Abder.Touati@fsb.rnu.tn